\documentclass[11pt,a4paper]{article}
\usepackage[margin=1in]{geometry}
\usepackage{latexsym, amsfonts, amsmath,tikz}
\usepgflibrary{arrows}
\newcommand{\be}{\begin{equation}}
\newcommand{\ee}{\end{equation}}
\newcommand{\bea}{\begin{eqnarray}}
\newcommand{\eea}{\end{eqnarray}}
\newcommand{\bean}{\begin{eqnarray*}} 
\newcommand{\eean}{\end{eqnarray*}}

\newcommand{\brray}{\begin{array}}
\newcommand{\erray}{\end{array}}
\newcommand{\ben}{\begin{equation}{nonumber}}
\newcommand{\een}{\end{equation}{nonumber}}

\newtheorem{dfn}{Definition}[section]
\newtheorem{thm}[dfn]{Theorem}
\newtheorem{lmma}[dfn]{Lemma}
\newtheorem{ppsn}[dfn]{Proposition}
\newtheorem{crlre}[dfn]{Corollary}
\newtheorem{xmpl}[dfn]{Example}
\newtheorem{rmrk}[dfn]{Remark}
\newcommand{\bdfn}{\begin{dfn}}
\newcommand{\bthm}{\begin{thm}}
\newcommand{\blmma}{\begin{lmma}}
\newcommand{\bppsn}{\begin{ppsn}}
\newcommand{\bcrlre}{\begin{crlre}}
\newcommand{\bxmpl}{\begin{xmpl}}
\newcommand{\brmrk}{\begin{rmrk}}
\newcommand{\edfn}{\end{dfn}}
\newcommand{\ethm}{\end{thm}}
\newcommand{\elmma}{\end{lmma}}
\newcommand{\eppsn}{\end{ppsn}}
\newcommand{\ecrlre}{\end{crlre}}
\newcommand{\exmpl}{\end{xmpl}}
\newcommand{\ermrk}{\end{rmrk}}

\newcommand{\cla}{{\cal A}}
\newcommand{\clb}{{\cal B}}
\newcommand{\clc}{{\cal C}}
\newcommand{\cld}{{\cal D}}

\newcommand{\clf}{{\cal F}}

\newcommand{\cli}{{\cal I}}
\newcommand{\clj}{{\cal J}}

\newcommand{\cln}{{\cal N}}
\newcommand{\clo}{{\cal O}}

\newcommand{\clq}{{\cal Q}}

\newcommand{\clt}{{\cal T}}

\newcommand{\clv}{{\cal V}}
\newcommand{\clw}{{\cal W}}

\def\a*{{\cal A}_{h,*}}
\def\B{{\cal B}(h)}
\def\B1{{\cal B}_1(h)}
\def\b{{\cal B}^{\rm s.a.}(h)}
\def\b1{{\cal B}^{\rm s.a.}_1(h)}

\newcommand{\ot}{\otimes}

\newcommand{\raro}{\rightarrow}

\def \qed {$\Box$}

\newcommand{\midarrow}{\tikz \draw[-triangle 90] (0,0) -- +(.1,0);}

\def\a*{{\cal A}_{h,*}}
\def\B{{\cal B}(h)}
\def\B1{{\cal B}_1(h)}
\def\b{{\cal B}^{\rm s.a.}(h)}
\def\b1{{\cal B}^{\rm s.a.}_1(h)}

\begin{document}
\begin{center}
{\Large{\bf Quantum symmetry of graph $C^{\ast}$-algebras at critical inverse temperature }}\\ 
{\large {\bf Soumalya Joardar \footnote{Acknowledges support from DST INSPIRE faculty grant} and \bf Arnab Mandal \footnote{Acknowledges support from SERB}}}\\ 
\end{center}
\begin{abstract}
We give a notion of quantum automorphism group of graph $C^{\ast}$-algebra without sink at critical inverse temperature. This is defined to be the universal object of a category of CQG's having a linear action in the sense of \cite{sou_arn} and preserving the KMS state at critical inverse temperature. We show that this category for a certain KMS state at critical inverse temperature coincides with the category introduced in \cite{sou_arn} for a class of graphs. We also introduce an orthogonal filtration on Cuntz algebra with respect to the unique KMS state and show that the category of CQG's preserving the orthogonal filtration coincides with the category introduced in this paper.
 
\end{abstract}

 \section{Introduction}
 Since S. Wang introduced the concept of quantum automorphism group of finite spaces in his pioneering paper \cite{Wang}, study of quantum automorphism groups of various noncommutative topological and geometric structures has come a long way. The idea was to formulate the classical group symmetry problem in categorical language i.e. to realize the automorphism group of the underlying space as a universal object in a certain category of compact groups and then letting the objects in the category be replaced by compact quantum groups and the group homomorphisms by CQG morphisms. Then the main challenge was to prove the existence of universal object in the larger category. As it turned out even for an $n$-point space where the underlying $C^{\ast}$-algebra is a finite dimensional $C^{\ast}$-algebra, the quantum automorphism is significantly larger than the classical symmetry group $S_{n}$ and in fact as a $C^{\ast}$-algebra the quantum automorphism group is infinite dimensional for $n\geq 4$. In that paper he also defined quantum automorphism group for finite dimensional $C^{\ast}$-algebras. And in deed for such $C^{\ast}$-algebras, the universal object fails to exist in general. S. Wang remedied this problem by his pioneering idea of restricting the category of $C^{\ast}$-actions of CQG's by demanding that such an action should also preserve some suitable linear functional.  Later  the concept of quantum automorphism groups was extended for finite graphs by T. Banica,  J. Bichon and others (see \cite{Ban}, \cite{Bichon} and references therein). More recently the notion of quantum automorphism groups was extended in the set up of NCG {\it ala} Connes (see \cite{connes}) by D. Goswami and others (see \cite{Debashish},\cite{Laplace} and references therein) and such notion was in a sense generalized by A.Skalski and T. Banica in \cite{Ortho}, where they have introduced the notion of orthogonal filtration on a $C^{\ast}$-algebra equipped with a faithful state.\\
 \indent Among the noncommutative topological objects, graph $C^{\ast}$-algebras are fairly well studied. But surprisingly no one really studied quantum automorphism groups of graph $C^{\ast}$-algebras until very recently S. Schmidt and M. Weber gave a formulation of quantum automorphism groups of graph $C^{\ast}$-algebras in \cite{Web}. Their framework, however, was algebraic. Following their work the authors of this paper gave a more analytic formulation of the same in \cite{sou_arn}. The formulations differ significantly in the sense that the quantum automorphism group of the underlying graphs turn out to be a subobject of the quantum automorphoism groups in \cite{sou_arn}, whereas in \cite{Web} the quantum automorphism group turns out to be the same as that of the underlying graph in the sense of \cite{Ban}. In their formulation in \cite{sou_arn}, the authors, motivated by the idea of S. Wang, defined a category of CQG's acting linearly on graph $C^{\ast}$-algebras and preserving a certain linear functional on a certain vector space. In this paper  we seek for a state on the whole graph $C^{\ast}$-algebras and define the quantum automorphism group to preserve the state on the whole $C^{\ast}$-algebra.\\
 \indent  The study of KMS states is quite old and dates back to 1980 where KMS state on Cuntz algebra was studied in \cite{ped}. Then for general graph $C^{\ast}$-algebras the canonical gauge invariant KMS states were studied by Laca, Watatani, Enomoto and others  (see \cite{Laca},\cite{watani} for example). In this paper we concentrate on graph $C^{\ast}$-algebras without sink and KMS states on such graph $C^{\ast}$-algebras at inverse critical temperature. Under certain hypothesis on the vertex matrix of the underlying graphs, such KMS state exists. Then  we take our quantum automorphism group to be the universal object in the category of CQG's acting linearly and preserving such gauge invariant KMS state at critical inverse temperature. Then we show that for a large class of graphs and a certain KMS state at inverse critical temperature, the category considered in this paper and the category considered in \cite{sou_arn} are same, simplifying the computation of quantum automorphism groups to a large extent. Also in case of Cuntz algebra we introduce an orthogonal filtration with respect to the unique gauge invariant ${\rm KMS}_{{\rm log}n}$ state and show that the corresponding orthogonal filtration preserving category and our category are same.\\
\indent Now let us briefly discuss the organisation of the paper. The paper starts with a preliminary section where we recall some basic facts about compact quantum groups, quantum symmetry and KMS states on graph $C^{\ast}$-algebras. Mainly, KMS states at critical inverse temperature on graph $C^{\ast}$-algebras is discussed. While recalling quantum symmetry, we also recall the notion of quantum symmetry of $C^{\ast}$-algebras equipped with an orthogonal filtration with respect to a faithful state from \cite{Ortho}. Then in the 3rd section a formulation of quantum automorphism group of graph $C^{\ast}$-algebras at critical inverse temperature is given. This is defined to be the universal object of a certain category introduced in this paper. We specialize to a particular class of graphs and it is shown that in fact the new category introduced in this paper for a certain KMS state coincides with the category introduced in \cite{sou_arn}. This helps in simplifying computation of quantum automorphism groups of a large class of graph $C^{\ast}$-algebras. Then in the next section two examples of quantum automorphism groups are computed. In the case of a graph with $n$-disjoint leaves, for a choice of KMS state, the quantum automorphism group at critical inverse temperature turns out to be an easy quantum group discussed in \cite{easy}. Then in the last subsection, we concentrate on quantum symmetry of Cuntz algebra at well known critical inverse ${\rm KMS}_{{\rm log}n}$ state $\phi$. There an orthogonal filtration with respect to $\phi$ is introduced and it is shown that the quantum automorphism group preserving that particular orthogonal filtration coincides with the quantum automorphism group at critical inverse temperature in our sense which turns out to be the well known quantum group $U^{+}_{n}$.  
\section{Preliminaries}
\subsection{Notational convention} 
In this paper all the graphs are {\bf finite, directed}. All the $C^{\ast}$-algebras are unital. For a $C^{\ast}$-algebra $\cla$, $\cla^{\ast}$ will denote the space of all linear functionals on $\cla$. The tensor product $\ot$ between two $C^{\ast}$-algebras are injective tensor products. ${\rm id}$ will denote the identity map. For a set $A$, ${\rm Sp}(A)$ ($\overline{\rm Sp}(A)$) will denote the linear (closed linear) span of the elements of $A$. $\mathbb{N}$ will denote the set of natural numbers and $\mathbb{N}_{0}$ will denote the set $\mathbb{N}\cup\{0\}$. By a universal object in a category we shall always mean the initial object in the category.
 \subsection{Compact quantum groups and quantum automorphism groups}
\label{qaut}
In this subsection we recall the basics of compact quantum groups and their actions on $C^{\ast}$-algebras. The facts collected in this Subsection are well known and we refer the readers to \cite{Van}, \cite{Woro}, \cite{Wang} for details. 
\bdfn
A compact quantum group (CQG) is a pair $(\clq,\Delta)$ such that $\clq$ is a unital $C^{\ast}$-algebra and $\Delta:\clq\raro\clq\ot\clq$ is a unital $C^{\ast}$-homomorphism satisfying\\
(i) $({\rm id}\ot\Delta)\circ\Delta=(\Delta\ot{\rm id})\circ\Delta$.\\
(ii) {\rm Sp}$\{\Delta(\clq)(1\ot\clq)\}$ and {\rm Sp}$\{\Delta(\clq)(\clq\ot 1)\}$ are dense in $\clq\ot\clq$.
\edfn
Given a CQG $\clq$, there is a canonical dense Hopf $\ast$-algebra $\clq_{0}$ in $\clq$ on which an antipode $\kappa$ and counit $\epsilon$ are defined. Given two CQG's $(\clq_{1},\Delta_{1})$ and $(\clq_{2},\Delta_{2})$, a CQG morphism between them is a $C^{\ast}$-homomorphism $\pi:\clq_{1}\raro\clq_{2}$ such that $(\pi\ot\pi)\circ\Delta_{1}=\Delta_{2}\circ\pi$.
\vspace{0.1in}\\
{\it Examples}:\\
1. Let $U_{n}^{+}$ be the universal $C^{\ast}$-algebra generated by $n^2$ elements $\{q_{ij}\}_{i,j=1,\ldots,n}$ such that both the matrices $((q_{ij}))$ and $((q_{ij}))^t$ are unitaries. Then $U_{n}^{+}$ is a CQG with coproduct given on the generators by $\Delta_{0}(q_{ij})=\sum_{k=1}^{n}q_{ik}\ot q_{kj}$ (see \cite{Wangfree} where this CQG has been denoted by $A_{u}(n)$).
\vspace{0.1in}\\
2. We discuss a particular unitary easy compact quantum group. For details on unitary easy compact quantum groups the reader might consult \cite{easy}. The unitary easy CQG $H_{n}^{\infty +}$ is defined to be the universal $C^{\ast}$-algebra generated by $\{q_{ij}:i,j=1,...,n\}$ such that\\
(a) the matrices $((q_{ij}))$ and $((q_{ij}^{\ast}))$ are unitaries,\\
(b) each $q_{ij}$'s are normal partial isometries.\\
The coproduct $\Delta$ is given on the generators by $\Delta(q_{ij})=\sum_{k=1}^{n}q_{ik}\ot q_{kj}$.\\
We record the following Lemma for future purpose.
\blmma
\label{H}
For $q_{ij}\in H_{n}^{\infty +}$, $q_{ki}q_{kj}=0$ for $i\neq j$ and $k=1,...,n$.
\elmma
{\it Proof}:\\
Using the fact that  for fixed $k$, $q_{ik}q_{ik}^{\ast}$'s are mutually orthogonal projections for $i=1,...,n$, one can show that $q_{ik}q_{jk}=0$ for $i\neq j$ (see discussion following the definition of quantum reflection group in \cite{easy}). Applying $\kappa$, we prove the Lemma.\qed
\vspace{0.1in}\\
3. For the following construction see \cite{Soltan}. Let $(\mathcal{Q},\Delta)$ be a CQG with a CQG-automorphism $\theta$ such that $\theta^2={ \rm Id}$. The doubling of this CQG, 
say $(\mathcal{D}_{\theta}(\mathcal{Q}),\tilde{\Delta})$, is  given by  $\mathcal{D}_{\theta}(\mathcal{Q}) :=\mathcal{Q}\oplus \mathcal{Q}$ (direct sum as a $C^{\ast}$-algebra),
and the coproduct is defined by the following
\begin{eqnarray*}
	&& \tilde{\Delta} \circ \xi= (\xi \ot \xi + \eta \ot [\eta \circ \theta])\circ  \Delta\\
	&& \tilde{\Delta} \circ \eta= (\xi \ot \eta + \eta \ot [\xi \circ \theta])\circ \Delta,
\end{eqnarray*}
where we have denoted  the injections of $\mathcal{Q}$ onto 
the first and second coordinate in $\mathcal{D}_{\theta}(\mathcal{Q})$ by $\xi$ and $\eta$ respectively, i.e. 
$\xi(a)=(a,0), \  \eta(a)= (0,a), \ (a \in \mathcal{Q}).$
\bdfn
Given a (unital) $C^{\ast}$-algebra $\cla$, a CQG $(\clq,\Delta)$ is said to act faithfully on $\cla$ if there is a unital $C^{\ast}$-homomorphism $\alpha:\cla\raro\cla\ot\clq$ satisfying\\
(i) $(\alpha\ot {\rm id})\circ\alpha=({\rm id}\ot \Delta)\circ\alpha$.\\
(ii) {\rm Sp}$\{\alpha(\cla)(1\ot\clq)\}$ is dense in $\cla\ot\clq$.\\
(iii) The $\ast$-algebra generated by the set  $\{(\omega\ot{\rm id})\alpha(\cla): \omega\in\cla^{\ast}\}$ is norm-dense in $\clq$.
\edfn
For a faithful action of a CQG $(\clq,\Delta)$ on a unital $C^{\ast}$-algebra $\cla$, there is a norm dense $\ast$-subalgebra $\cla_{0}$ of $\clc$ such that the canonical Hopf-algebra $\clq_{0}$ coacts on $\cla_{0}$.
\bdfn
(Def 2.1 of \cite{Bichon})
Given a unital $C^{\ast}$-algebra $\cla$, quantum automorphism group of $\cla$ is a CQG $(\clq,\Delta)$ acting faithfully on $\cla$ satisfying the following universal property:\\
\indent If $(\clb,\Delta_{\clb})$ is any CQG acting faithfully on $\cla$, there is a surjective CQG morphism $\pi:\clq\raro\clb$ such that $({\rm id}\ot \pi)\circ\alpha=\beta$, where $\beta:\cla\raro\clc\ot\clb$ is the corresponding action of $(\clb,\Delta_{\clb})$ on $\cla$ and $\alpha$ is the action of $(\clq,\Delta)$ on $\cla$.
\edfn
\brmrk
In general universal object might fail to exist in the above category. For the existence of universal object one generally restricts the category. In addition to $\alpha$ being a faithful action it is common to add some kind of `volume' preserving condition i.e. one assumes some linear functional $\tau$ on $\clc$ such that $(\tau\ot {\rm id})\circ\alpha(a)=\tau(a).1$ for all $a$. in some suitable subspace of $\clc$ (see \cite{Wang}, \cite{Debashish}).
\ermrk
{\it Example}:\\
1. If we take a space of $n$ points $X_{n}$ then the quantum automorphism group of the $C^{\ast}$-algebra $C(X_{n})$ is given by the CQG (denoted by $S_{n}^{+}$) which as a $C^{\ast}$-algebra is the universal $C^{\ast}$ algebra generated by $\{q_{ij}\}_{i,j=1,\ldots,n}$ satisfying the following relations (see Theorem 3.1 of \cite{Wang}):
\begin{displaymath}
 q_{ij}^{2}=q_{ij},q_{ij}^{\ast}=q_{ij},\sum_{k=1}^{n}q_{ik}=\sum_{k=1}^{n}q_{ki}=1, i,j=1,...,n.
\end{displaymath}
The coproduct on the generators is given by $\Delta(q_{ij})=\sum_{k=1}^{n}q_{ik}\ot q_{kj}$.\vspace{0.1in}\\
2. If we take the $C^{\ast}$-algebra $M_{n}(\mathbb{C})$, then as remarked earlier if we take the category of CQG's only having faithful $C^{\ast}$-action on $M_{n}(\mathbb{C})$, then universal object fails to exist in the category. But in addition if we assume any object in the category also has to preserve a linear functional $\phi$ on $M_{n}(\mathbb{C})$, then existence of the universal object can be shown (see \cite{Wang}). 
 \subsection{Orthogonal filtration preserving quantum symmetry}
\label{orthogonalfiltration}
 In this Subsection, we recall the formulation of quantum symmetry of a $C^{\ast}$-algebra equipped with an orthogonal filtration with respect to a state on the $C^{\ast}$-algebra from \cite{Ortho}.
 \bdfn
\label{filtration}
 Let $\cla$ be a unital $C^{\ast}$-algebra equipped with a faithful state $\phi$ and a family $\{V_{i}\}_{i\in\cli}$ of finite dimensional subspaces of $\cla$ (with the index set $\cli$ containing a distinguished element $0$) satisfying the following conditions:\\
 (i) $V_{0}=\mathbb{C}1_{\cla}$;\\
 (ii) for all $i,j\in\cli$, $i\neq j$, $a\in V_{i},b\in V_{j}$, we have $\phi(a^{\ast}b)=0$;\\
 (iii) the set ${\rm Sp}(\cup_{i\in\cli}V_{i})$ is dense in $\cla$.\\
 If the above conditions are satisfied then we say that the pair $(\phi,(V_{i})_{i\in\cli})$ defines an orthogonal filtration on $\cla$.
 \edfn
 \bdfn
 Let $\cla$ be a unital $C^{\ast}$-algebra with an orthogonal filtration $(\phi,(V_{i})_{i\in \cli})$. We say that a CQG $(\clq,\Delta)$ acts on $\cla$ in a filtration preserving way if there is $C^{\ast}$-action $\alpha$ of $\clq$ on $\cla$ such that 
 \begin{displaymath}
 \alpha(V_{i})\in V_{i}\ot \clq,  i\in\cli.
 \end{displaymath} We shall then write $(\alpha,(\clq,\Delta))\in\clc_{\cla,\mathcal{V}}$.
 \edfn
 \brmrk
\label{statepreserve}
 If $(\alpha,(\clq,\Delta))\in\clc_{\cla,\mathcal{V}}$, then $(\phi\ot{\rm id})\circ\alpha(a)=\phi(a)1_{\cla}$ for all $a\in\cla$. It follows from the simple observation that for $a\in {\rm Sp}(\cup_{i\in\cli-\{0\}}V_{i})$, $\phi(a)=0$.
 \ermrk
 We have the following (Theorem 2.7 of \cite{Ortho})
 \bthm
 Let $(\cla,\phi,(V_{i})_{i\in\cli})$ be a $C^{\ast}$-algebra with orthogonal filtration. Then there exists universal object in the category $\clc_{\cla,\mathcal{V}}$.
 \ethm
 The universal object in the category $\clc_{\cla,\mathcal{V}}$ is said to be the quantum automorphism group of $(\cla,\phi,(V_{i})_{i\in\cli})$.
\subsection{KMS states on graph $C^{\ast}$-algebra without sink}
A {\bf finite} directed graph is a collection of finitely many edges and vertices. If we denote the edge set of a graph $\Gamma$ by $E=(e_{1},\ldots,e_{n})$ and set of vertices of $\Gamma$ by $V=(v_{1},\ldots,v_{m})$ then recall the maps $s,t:E\raro V$ and the vertex matrix $D$ is an $m\times m$ matrix whose $ij$-th entry is $k$ if there are $k$-number of edges from $v_{i}$ to $v_{j}$. We also denote the space of paths by $E^{\ast}$ (see \cite{Laca}).
\bdfn
$\Gamma$ is said to be without sink if the map $s:E\raro V$ is surjective. Furthermore $\Gamma$ is said to be without multiple edges if the adjacency matrix $D$ has entries either $1$ or $0$. 
\edfn
\brmrk
Note that the graph $C^{\ast}$-algebra corresponding to a graph without sink is a Cuntz-Krieger algebra. Reader might see \cite{Cuntz-Krieger} for more details on Cuntz-Krieger algebra.
\ermrk
Now we recall some basic facts about graph $C^{\ast}$-algebras. Reader might consult \cite{Pask} for details on graph $C^{\ast}$-algebras. Let $\Gamma=\{E=(e_{1},...,e_{n}),V=(v_{1},...,v_{m})\}$ be a finite graph without sink. A graph is said to be {\bf connected} if every vertex is a source or a target of an edge. In other words for all $v_{l}\in V$, $s^{-1}(v_{l})$ or $t^{-1}(v_{l})$ is non empty. In this paper all the graphs are {\bf finite, connected} and without sink. We assign partial isometries $S_{i}$'s to edges $e_{i}$ for all $i=1,...,n$ and projections $p_{v_{i}}$ to the vertices $v_{i}$ for all $i=1,...,m$.
\bdfn
\label{Graph}
The graph $C^{\ast}$-algebra $C^{\ast}(\Gamma)$ is defined as the universal $C^{\ast}$-algebra generated by partial isometries $\{S_{i}\}_{i=1,\ldots,n}$ and mutually orthogonal projections $\{p_{v_{k}}\}_{k=1,\ldots,m}$ satisfying the following relations:
\begin{displaymath}
  S_{i}^{\ast}S_{i}=p_{t(e_{i})}, \sum_{s(e_{j})=v_{l}}S_{j}S_{j}^{\ast}=p_{v_{l}}.
	\end{displaymath}
\edfn
In a graph $C^{\ast}$-algebra $C^{\ast}(\Gamma)$, we have the following (see Subsection 2.1 of \cite{Pask}):\\
1. $\sum_{k=1}^{m}p_{v_{k}}=1$ and $S_{i}^{\ast}S_{j}=0$ for $i\neq j$.
\vspace{0.05in}\\
2. $S_{\mu}=S_{1}S_{2}\ldots S_{l}$ is non zero if and only if $\mu=e_{1}e_{2}\ldots e_{l}$ is a path i.e. $t(e_{i})=s(e_{i+1})$ for $i=1,\ldots,(l-1)$.
\vspace{0.05in}\\
3. $C^{\ast}(\Gamma)={\overline{\rm Sp}}\{S_{\mu}S_{\nu}^{\ast}:t(\mu)=t(\nu)\}$. We denote ${\rm Sp}\{S_{\mu}S_{\nu}^{\ast}:t(\mu)=t(\nu)\}$ by $\cla_{0}$. It can be shown that $\cla_{0}$ is a $\ast$-subalgebra of $C^{\ast}(\Gamma)$.\\
Now we shall briefly discuss KMS-states on graph $C^{\ast}$-algebras. For that we recall Toeplitz algebra $\clt C^{\ast}(\Gamma)$. Readers are referred to \cite{Laca} for details. Our convention though is opposite to that of \cite{Laca} in the sense that we interchange source projections and target projections. Also we shall modify the results of \cite{Laca} according to our need. Suppose that $\Gamma$ is a directed graph as before. A Toeplitz-Cuntz-Krieger $\Gamma$ family consists of mutually orthogonal projections $\{p_{v_{i}}:v_{i}\in V\}$ and partial isometries $\{S_{i}:e_{i}\in E\}$ such that $\{S_{i}^{\ast}S_{i}=p_{t(e_{i})}\}$ and
\begin{displaymath}
p_{v_{l}}\geq \sum_{s(e_{i})=v_{l}}S_{i}S_{i}^{\ast}.
\end{displaymath}
Toeplitz algebra $\clt C^{\ast}(\Gamma)$ is defined to be the universal $C^{\ast}$-algebra generated by the Toeplitz-Cuntz-Krieger $\Gamma$ family. It is clear from the definition that $C^{\ast}(\Gamma)$ is the quotient of $\clt C^{\ast}(\Gamma)$ by the ideal $\clj$ generated by
\begin{displaymath}
  P:=\{p_{v_{l}}-\sum_{s(e_{i})=v_{l}}S_{i}S_{i}^{\ast}\}.
	\end{displaymath}
	The standard arguments give $\clt C^{\ast}(\Gamma)=\overline{\rm Sp} \{S_{\mu}S_{\nu}^{\ast}:t(\mu)=t(\nu)\}$. $\clt C^{\ast}(\Gamma)$ admits the usual gauge action $\gamma$ of $\mathbb{T}$ which descends to the usual gauge action on $C^{\ast}(\Gamma)$ given on the generators by $\gamma_{z}(S_{\mu}S_{\nu}^{\ast})=z^{(|\mu|-|\nu|)}S_{\mu}S_{\nu}^{\ast}$. Consequently it has a dynamics $\alpha:\mathbb{R}\raro {\rm Aut} \ C^{\ast}(\Gamma)$ which is lifted from $\gamma$ via the map $t\raro e^{it}$. We recall the following from \cite{Laca} (Proposition 2.1)
	\bppsn
	\label{exist_KMS}
	Let $\Gamma$ be a finite, directed, connected graph without sink and $\gamma:\mathbb{T}\raro {\rm Aut} \ \clt C^{\ast}(\Gamma)$ be the gauge action with the corresponding dynamics $\alpha:\mathbb{R}\raro {\rm Aut} \ \clt C^{\ast}(\Gamma)$. Let $\beta\in\mathbb{R}$.\\
	(a) A state $\phi$ is a ${\rm KMS}_{\beta}$ state of $(\clt C^{\ast}(\Gamma),\alpha)$ if and only if
	\begin{displaymath}
	\phi(S_{\mu}S_{\nu}^{\ast})=\delta_{\mu,\nu}e^{-\beta|\mu|}\phi(p_{t(\mu)}).
	\end{displaymath}
	(b) Suppose that $\phi$ is a ${\rm KMS}_{\beta}$ state of $(\clt C^{\ast}(\Gamma),\alpha)$, and define $\cln^{\phi}=(\cln^{\phi}_{i})\in[0,\infty)^{m}$ by $\cln^{\phi}_{i}=\phi(p_{v_{i}})$. Then $\cln^{\phi}$ is a probability measure on $V$ satisfying the subinvariance condition $D\cln^{\phi}\leq e^{\beta}\cln^{\phi}$.\\
	(c) A ${\rm KMS}_{\beta}$ state factors through $C^{\ast}(\Gamma)$ if and only if $(D\cln^{\phi})_{i}=e^{\beta}\cln^{\phi}_{i}$ for all $i=1,\ldots, m$ i.e. $\cln^{\phi}$ is an eigen vector of $D$ with eigen value $e^{\beta}$.
	\eppsn
	{\it Proof}:\\
	We only prove part (c) as we are concerned about graph $C^{\ast}$-algebras in this paper. Proofs of all the statements can be found in \cite{Laca}. Since we have assumed the graph is without sink, for any $v_{l}\in V$, $s^{-1}(v_{l})\neq\emptyset$. Any state on $\clt C^{\ast}(\Gamma)$ that vanishes on the ideal $\clj$, descends to $C^{\ast}(\Gamma)$. Now $\phi(p_{v_{l}}-\sum_{s(e_{i})=v_{l}}S_{i}S_{i}^{\ast})=\phi(p_{v_{l}})-e^{-\beta}\sum_{s(e_{i})=v_{l}}\phi(p_{t(e_{i})})$. That is the state descends to $C^{\ast}(\Gamma)$ if and only if $e^{\beta}\cln^{\phi}_{l}=\sum_{v_{k}\in V}D(l,k)\cln^{\phi}_{k}$.\qed 
\subsubsection{KMS state at critical inverse temperature}
\label{KMS}
In this subsection we collect a few results on existence of KMS states on graph $C^{\ast}$-algebras. For that we continue to assume $\Gamma$ to be a finite, connected graph without sink and with vertex matrix $D$. We denote the spectral radius of $D$ by $\rho(D)$. We start with the following proposition (Proposition 4.1 of \cite{Laca}) which gives a condition for the existence of ${\rm KMS}_{\beta}$ state at critical inverse temperature $\rho(D)$.
\bppsn
\label{exist1}
If $\cln$ is a probability measure on $V$ such that $D\cln\leq \rho(D)\cln$, then there is a ${\rm KMS}_{ln\rho(D)}$ state on $\clt C^{\ast}(\Gamma)$ such that
\begin{displaymath}
\phi(S_{\mu}S_{\nu}^{\ast})=\delta_{\mu,\nu}\rho(D)^{-|\mu|}\cln_{i},
\end{displaymath} 
where $t(\mu)=v_{i}$. The state factors through a state to $C^{\ast}(\Gamma)$ if and only if $D\cln=\rho(D)\cln$.
\eppsn 
The above proposition has the following (Corollary 4.2 of \cite{Laca})
\bcrlre
Let $\Gamma$ be a finite directed, connected graph without sink and with vertex matrix $D$. Then $(\clt C^{\ast}(\Gamma),\alpha)$ has a ${\rm KMS}_{ln\rho(D)}$ state.
\ecrlre
In the following we give a necessary and sufficient condition on the graph $\Gamma$ for the existence of ${\rm KMS}_{ln\rho(D)}$ state on $C^{\ast}(\Gamma)$. Before that we recall the fixed point subalgebra of $C^{\ast}(\Gamma)$ with respect to the canonical gauge action given on the generators by
\begin{displaymath}
\gamma_{z}(S_{i})=zS_{i}, i=1,...,n.
\end{displaymath}
It is well known that the fixed point subalgebra $C^{\ast}(\Gamma)^{\gamma}$ is an AF algebra with the finite dimensional $C^{\ast}$-subalgebras given by
\begin{displaymath}
\clf_{0}=\mathbb{C}1, \clf_{k}={\rm Sp}\{S_{\mu}S_{\nu}^{\ast}:|\mu|=|\nu|=k:k\geq 1\}.
\end{displaymath}
$\clf_{k}\subset \clf_{k+1}$ with the embedding given by $S_{\mu}S_{\nu}^{\ast}=S_{\mu}(\sum_{j=1}^{n}S_{j}S_{j}^{\ast})S_{\nu}^{\ast}$. Note that here we have used the fact that the graph has no sink so that $\sum_{j}S_{j}S_{j}^{\ast}=1$.
\bppsn
\label{exist_inverse}
$C^{\ast}(\Gamma)$ has a ${\rm KMS}_{ln\rho(D)}$ state (say $\phi$) which is faithful on $\clf_{k}$ for all $k$ if and only if $\rho(D)$ is an eigen value of $D$ with eigen vector $(w_{1},...,w_{m})$ such that $w_{i}>0$ for all $i=\{1,...,m\}$. If we normalize $(w_{i})$ so that $\sum w_{i}=1$, then the KMS state $\phi$ satisfies
\begin{displaymath}
\phi(S_{\mu}S_{\nu}^{\ast})=\delta_{\mu,\nu}\rho(D)^{-|\mu|}w_{i},
\end{displaymath}
where $t(\mu)=v_{i}$ and for $i=1,...,m$.
\eppsn
{\it Proof}:\\
Let $\rho(D)$ be an eigen value of $D$ with a normalized eigen vector $(w_{1},...,w_{m})$ such that $\sum w_{i}=1$ and $w_{i}>0$ for all $i=1,...,m$. Then the vector $\cln^{\phi}$ with $\cln^{\phi}_{i}=w_{i}$ is a probability measure on the vertex set such that $D\cln^{\phi}=\rho(D)\cln^{\phi}$. Then by Proposition \ref{exist1}, $C^{\ast}(\Gamma)$ has a ${\rm KMS}_{ln\rho(D)}$ state $\phi$ given by $\phi(S_{\mu}S_{\nu}^{\ast})=\delta_{\mu,\nu}\rho(D)^{-|\mu|}\cln^{\phi}_{i}=\delta_{\mu,\nu}\rho(D)^{-|\mu|}w_{i}$. From this formula it is clear that $\phi$ is faithful on $\clf_{k}$ for all $k$. The converse follows from $(c)$ of Proposition \ref{exist_KMS}. The faithfulness ensures that all the entries of the eigen vector is positive.\qed\vspace{0.1in}\\
{\it Examples}:\\
1. {\bf Strongly connected graphs}: A graph is said to be strongly connected if $vE^{\ast}w$ is non empty for all $v,w\in V$. By Theorem 1.5 of \cite{Seneta}, for a strongly connected graph with vertex matrix $D$, $\rho(D)$ is an eigen value of $D$ for which there is an eigen vector with all positive entries called the Perron-Frobenius eigen vector of $D$. Examples of such graphs include graph of Cuntz algebra, polygon, complete graph in the sense of \cite{sou_arn}. Note that in case of strongly connected graphs, in fact there is a unique ${\rm KMS}$ state (see Theorem 4.3 of \cite{Laca}).\vspace{0.1in} \\
2. {\bf Disjoint n- leaves}: This is a graph with $n$-vertices such that each vertex has a loop. The corresponding vertex matrix is given by the identity matrix and hence it has ${\rm KMS}_{0}$ state corresponding to every probability measure given on the vertex set. This is because the eigen space corresponding to the eigen value $1$ is whole of $\mathbb{R}^{n}$. 
\section{Quantum symmetry}
\bdfn
Let $\Gamma$ be a finite, directed, connected graph without sink and with vertex matrix $D$ such that $\Gamma$ has a ${\rm KMS}_{ln\rho(D)}$ state $\phi$. We define the category $\clc_{\phi}^{\Gamma}:=(\alpha,(\clq,\Delta))$ where $(\clq,\Delta)$ is a CQG having linear action on $C^{\ast}(\Gamma)$ such that it preserves $\phi$ in the following sense:
\begin{eqnarray*}
&&\alpha(S_{i})=\sum_{j=1}^{n}S_{j}\ot q_{ji}, q_{ij}\in \clq,\\
&&(\phi\ot {\rm id})\circ\alpha(a)=\phi(a).1, a\in C^{\ast}(\Gamma).
\end{eqnarray*}
Morphism between two objects $((\clq_{1},\Delta_{1}),\alpha_{1})$ and $((\clq_{2},\Delta_{2}),\alpha_{2})$ in the category is given by a CQG morphism $\Phi:\clq_{1}\raro\clq_{2}$ such that $({\rm id}\ot\Phi)\circ\alpha_{1}=\alpha_{2}$.
\edfn
We call any CQG action which preserves ${\rm Sp}\{S_{i}:i=1,...,n\}$ a linear action. Adapting the arguments of Lemma 5.4.2 of \cite{thesis}, we can prove the following
\bppsn
The category $\clc_{\phi}^{\Gamma}$ has a universal object.
\eppsn
Now we recall the category introduced in \cite{sou_arn}.
\bdfn (definition 4.7 of \cite{sou_arn})
Let $\Gamma$ be a finite, directed graph without sink. The category $\clc_{\tau}^{\rm Lin}$ has objects $(\alpha,(\clq,\Delta))$ where $\alpha$ is a linear action on $C^{\ast}(\Gamma)$ such that it preserves $\tau$, where $\tau$ is a linear functional on a vector space $\clv_{2}^{+}:={\rm Sp}\{S_{i}S_{j}^{\ast}:i,j=1,...,n\}$ given by $\tau(S_{i}S_{j}^{\ast})=\delta_{ij}$. The morphisms are the usual ones i.e. CQG morphisms intertwining the actions.  
\edfn
 \blmma
\label{nopath}
Let $\Gamma$ be a finite, directed graph such that $(\alpha,(\clq,\Delta))$ is an object in the category $\clc^{\rm Lin}_{\tau}$ with $\alpha(S_{i})=\sum_{j}S_{j}\ot q_{ji}$. For $S_{\mu_{1}}S_{\mu_{2}}...S_{\mu_{p}}\neq 0$ and  $S_{i_{1}}...S_{i_{p}}=0$, $q_{i_{1}\mu_{1}}...q_{i_{p}\mu_{p}}=0$.
\elmma
{\it Proof}:\\
Let $(i_{1},...,i_{p})$ be such that $S_{i_{1}}...S_{i_{p}}=0$. Then $\alpha(S_{i_{1}}...S_{i_{p}})=0$. On the other hand, 
\begin{displaymath}
\alpha(S_{i_{1}}...S_{i_{p}})=\sum S_{\gamma_{1}}...S_{\gamma_{p}}\ot q_{\gamma_{1}i_{1}}...q_{\gamma_{p}i_{p}}.
\end{displaymath}
Since $S_{\mu_{1}}S_{\mu_{2}}...S_{\mu_{p}}\neq 0$, we get $q_{\mu_{1}i_{1}}...q_{\mu_{p}i_{p}}=0$. Now recall from Proposition 4.8 of \cite{sou_arn} that the matrix $Q^{-1}=((q_{ij}))^{-1}$ is given by $Q^{-1}=(F^{\Gamma})^{-1}Q^{\ast}F^{\Gamma}$ for a diagonal matrix $F^{\Gamma}$. Then it follows that $\kappa(q_{ij})=c_{ij}q_{ji}^{\ast}$ for some non zero scalars $c_{ij}$ for all $i,j=1,...,n$. Now Applying $\kappa$ to the equation $q_{\mu_{1}i_{1}}...q_{\mu_{p}i_{p}}=0$, we conclude that $q_{i_{1}\mu_{1}}...q_{i_{p}\mu_{p}}=0$.\qed\vspace{0.1in}\\
Now from the following Lemma onwards we stick to a particular class of graphs whose vertex matrix $D$ is such that the row sums are all equal to $\rho(D)$.
\bthm
\label{newold}
As before let $\Gamma=\{E=(e_{1},...,e_{n}),V=(v_{1},...,v_{m})\}$ be a finite, directed, connected graph without sink and with vertex matrix $D$. If the row sums of $D$ are all equal to $\rho(D)$, then it has a ${\rm KMS}_{ln\rho(D)}$ state $\phi$ such that 
\begin{displaymath}
\phi(S_{\mu}S_{\nu}^{\ast})=\delta_{\mu,\nu}\rho(D)^{-|\mu|}\frac{1}{m}.
\end{displaymath}
In that case the categories $\clc^{\rm Lin}_{\tau}$ and $\clc^{\Gamma}_{\phi}$ are isomorphic.
\ethm
{\it Proof}:\\
Note that the assumption on the graph ensures that $\rho(D)$ is an eigenvalue of the matrix $D$ with one dimensional eigen subspace spanned by the vector $(1,...,1)$. Then we have a probability measure $\cln^{\phi}$ on the vertex set coming from the eigen subspace mentioned above given by $\cln^{\phi}_{i}=\frac{1}{m}$. Hence by Proposition \ref{exist_inverse} we have a ${\rm KMS}_{ln\rho(D)}$ state $\phi$ satisfying
\begin{displaymath}
\phi(S_{\mu}S_{\nu}^{\ast})=\delta_{\mu,\nu}\rho(D)^{-|\mu|}\frac{1}{m}.
\end{displaymath} 
Now let $(\alpha,(\clq,\Delta))$ be an object in the category $\clc^{\Gamma}_{\phi}$. Then it is easy to see that $\sum_{k=1}^{n}q_{ki}q_{kj}^{\ast}=\delta_{ij}.$ Recall the linear functional $\tau$ on $\clv^{+}_{2}$. Since the action is already linear, to show that $(\alpha,(\clq,\Delta))$ belongs to the category $\clc^{\rm Lin}_{\tau}$, it suffices to show that $(\tau\ot{\rm id})\circ\alpha(a)=\tau(a)1$ for all $a\in\clv^{+}_{2}$. To that end we consider $\alpha(S_{i}S_{j}^{\ast})=\sum S_{k}S_{l}^{\ast}\ot q_{ki}q_{lj}^{\ast}$. Using the fact that $\tau(S_{k}S_{l}^{\ast})=\delta_{kl}$, we get $(\tau\ot{\rm id})\circ\alpha(S_{i}S_{j}^{\ast})=\sum_{k}q_{ki}q_{kj}^{\ast}=\delta_{ij}=\tau(S_{i}S_{j}^{\ast})1$, proving that $(\alpha,(\clq,\Delta))$ belongs to the category $\clc^{\rm Lin}_{\tau}$. Conversely let $(\alpha,(\clq,\Delta))$ belong to the category $\clc^{\rm Lin}_{\tau}$. Then it follows that $\sum_{k}q_{ki}q_{kj}^{\ast}=\delta_{ij}$. It is easy to see that it suffices proving $(\phi\ot{\rm id})\circ\alpha(a)=\phi(a)1$ for all $a\in C^{\ast}(\Gamma)$. By linearity and continuity of the maps $\phi$ and $\alpha$ it is enough to show that 
\begin{displaymath}(\phi\ot{\rm id})\circ\alpha(S_{\mu}S_{\nu}^{\ast})=\phi(S_{\mu}S_{\nu}^{\ast}).1 \ \forall \ \mu,\nu\in E^{\ast}.
\end{displaymath}
 To that end, let $\mu=(\mu_{1},...,\mu_{p}),\nu=(\nu_{1},...,\nu_{r})$. Then if $p\neq r$, by definition of $\phi$, $\phi(S_{\mu}S_{\nu}^{\ast})=0$. Also $\alpha(S_{\mu}S_{\nu}^{\ast})=\sum S_{i_{1}}...S_{i_{p}}S_{j_{r}}^{\ast}...S_{j_{1}}^{\ast}\ot q_{i_{1}\mu_{1}}...q_{i_{p}\mu_{p}}q_{j_{r}\nu_{r}}^{\ast}...q_{j_{1}\nu_{1}}^{\ast}$. Again from definition of $\phi$, applying $(\phi\ot{\rm id}) $ on the R.H.S, we see that $(\phi\ot {\rm id})\circ\alpha(S_{\mu}S_{\nu}^{\ast})=0=\phi(S_{\mu}S_{\nu}^{\ast})1$. Now let $p=r$. Then $\phi(S_{i_{1}}...S_{i_{p}}S_{j_{p}}^{\ast}...S_{j_{1}}^{\ast})=\rho(D)^{-p}.\frac{1}{m}$ if $(i_{1},...,i_{p})=(j_{1},...,j_{p})$ and is equal to zero otherwise. We have
\begin{eqnarray*}
&&(\phi\ot{\rm id})\circ\alpha(S_{\mu}S_{\nu}^{\ast})\\
&=&(\phi\ot{\rm id}) (\sum S_{i_{1}}...S_{i_{p}}S_{j_{p}}^{\ast}...S_{j_{1}}^{\ast}\ot q_{i_{1}\mu_{1}}...q_{i_{p}\mu_{p}}q_{j_{p}\nu_{p}}^{\ast}...q_{j_{1}\nu_{1}}^{\ast})
\end{eqnarray*} 
Using the fact that if $\cld=\{(i_{1},...,i_{p}): S_{i_{1}}...S_{i_{p}}=0\}$, then $q_{i_{1}\mu_{1}}...q_{i_{p}\mu_{p}}=0$ for all $(i_{1},...,i_{p})\in\cld$ (Lemma \ref{nopath}), we get that the last summation equals to
\begin{displaymath}
\rho(D)^{-p}\frac{1}{m}\sum q_{i_{1}\mu_{1}}...q_{i_{p}\mu_{p}}q_{i_{p}\nu_{p}}^{\ast}...q_{i_{1}\nu_{1}}^{\ast}.
\end{displaymath} 
Hence for $(\mu_{1},...,\mu_{p})=(\nu_{1},...,\nu_{p})$, $(\phi\ot{\rm id})\circ\alpha(S_{\mu}S_{\mu}^{\ast})=\rho(D)^{-p}\frac{1}{m}\sum q_{i_{1}\mu_{1}}...q_{i_{p}\mu_{p}}q_{i_{p}\mu_{p}}^{\ast}...q_{i_{1}\mu_{1}}^{\ast}.$ Now using the fact that $\sum_{k}q_{ki}q_{ki}^{\ast}=1$ and interchanging the finite summations, we conclude that $(\phi\ot{\rm id})\circ\alpha(S_{\mu}S_{\mu}^{\ast})=\rho(D)^{-p}\frac{1}{m}=\phi (S_{\mu}S_{\mu}^{\ast})1$. Now let $(\mu_{1},...,\mu_{p})\neq(\nu_{1},...,\nu_{p})$ and without loss of generality, we assume that $\mu_{p}\neq \nu_{p}$. Then we get 
\begin{displaymath}
(\phi\ot{\rm id})\circ\alpha(S_{\mu}S_{\nu}^{\ast})=\rho(D)^{-p}\frac{1}{m}\sum q_{i_{1}\mu_{1}}...q_{i_{p}\mu_{p}}q_{i_{p}\nu_{p}}^{\ast}...q_{i_{1}\nu_{1}}^{\ast}.
\end{displaymath}
Now again using the fact $\sum_{i_{p}}q_{i_{p}\mu_{p}}q_{i_{p}\nu_{p}}^{\ast}=0$, for $\mu_{p}\neq\nu_{p}$, we get $(\phi\ot{\rm id})\circ\alpha(S_{\mu}S_{\nu}^{\ast})=0=\phi(S_{\mu}S_{\nu}^{\ast})1$. Hence $(\phi\ot{\rm id})\circ\alpha(S_{\mu}S_{\nu}^{\ast})=\phi(S_{\mu}S_{\nu}^{\ast})$ for all $\mu,\nu\in E^{\ast}$ which completes the proof of the theorem.\qed\vspace{0.1in}\\
Combining the above theorem with the results obtained in \cite{sou_arn}, we get the following for a finite, connected graph $\Gamma$ without sink and with vertex matrix $D$ such that the row sums are all equal to $\rho(D)$ (for definitions of $Q^{aut}_{Ban}(\Gamma)$ and $(C(S^{1})\star...\star C(S^{1}),\Delta)$ , see \cite{sou_arn} ),
\bcrlre
\label{old}
For the KMS state $\phi$ in Theorem \ref{newold},\\
(i) If we denote the universal object in the category $\clc^{\Gamma}_{\phi}$ by $Q^{\Gamma}_{\rm KMS}$, then $Q^{\Gamma}_{\rm KMS}\cong Q^{\rm Lin}_{\tau}$. We call $Q^{\Gamma}_{\rm KMS}$ to be the quantum automorphism group of the graph $\Gamma$ at inverse critical temperature.\\
(ii) $Q^{aut}_{Ban}(\Gamma)$ and $(C(S^{1})\star...\star C(S^{1}),\Delta)$ are subobjects of the category $\clc^{\Gamma}_{\phi}$.
\ecrlre
\section{Examples of quantum symmetry at inverse critical temperature}
1. {\bf Complete graph with 2 vertices}: Recall the complete graph $\Gamma$ with two vertices $\{0,1\}$ and two edges $e_{1},e_{2}$ from \cite{sou_arn} (see examle 2, page 11).
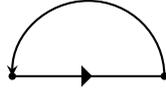
\begin{center} \begin{tikzpicture}
	\begin{scope}[thick, every node/.style={sloped,allow upside down}]
	\draw (0,7.5) -- node {\midarrow} (2,7.5);  
	\filldraw [black] (0,7.5) circle (1pt);
	\filldraw [black] (2,7.5) circle (1pt);
	\draw [black, -stealth ] (2,7.5) arc (0:180:1) node [right] {};
	\end{scope}
	\end{tikzpicture}
\end{center}
 If we denote the projections corresponding to the vertices by $p_{0},p_{1}$ and the partial isometries corresponding to the edges by $S_{1},S_{2}$, then $C^{\ast}(\Gamma)$ is generated by $S_{1},S_{2}$ satisfying the following relations:
\begin{eqnarray}
S_{1}^{\ast}S_{1}=S_{2}S_{2}^{\ast}=p_{1},S_{2}^{\ast}S_{2}=S_{1}S_{1}^{\ast}=p_{0}, 
\end{eqnarray}
where $p_{0},p_{1}$ are mutually orthogonal projections with $p_{0}+p_{1}=1$. It has the vertex matrix given by
$\begin{bmatrix}
	0 & 1\\
	1 & 0
	\end{bmatrix}$
. Hence the vertex matrix has its row sums equal to the spectral radius $1$. Since it is a strongly connected graph, it has a unique ${\rm KMS}_{0}$ state $\phi$, given by $\phi(S_{\mu}S_{\nu}^{\ast})=\delta_{\mu,\nu}\frac{1}{2}$. Combining Theorem 5.4 of \cite{sou_arn} and Theorem \ref{newold} of this paper we get the following result:
\bppsn
 $Q^{\Gamma}_{\rm KMS}\cong (\cld_{\theta}(C(S^{1})\star C(S^{1})),\tilde{\Delta})$
 \eppsn
    2. {\bf Disjoint n-leaves}: Let $\Gamma$ be a finite, connected, directed graph with vertex matrix ${\rm Id}_{n\times n}$. Then as discussed in Example 2 of Subsection \ref{KMS}, ${\rm KMS}_{0}$ states correspond to any probability measure given on the vertex set. $C^{\ast}(\Gamma)$ is generated by partial isometries $S_{i}$ satisfying
    \begin{eqnarray*}
    	S_{i}^{\ast}S_{i}=S_{i}S_{i}^{\ast}=p_{i}, p_{i}p_{j}=\delta_{ij}p_{i},p_{i}^{\ast}=p_{i},i,j=1,...,n.
    \end{eqnarray*} 
    For the following proposition recall the easy CQG $H_{n}^{\infty +}$ from Subsection \ref{qaut}.   
    \bppsn
    For the above graph $\Gamma$, the categories $\clc^{\rm Lin}_{\tau}$ and $\clc_{\phi}^{\Gamma}$ coincide if and only if the ${\rm KMS}_{0}$ state $\phi$ is given by
		\begin{displaymath}
		\phi(S_{\mu}S_{\nu}^{\ast})=\delta_{\mu,\nu}\frac{1}{n}.
		\end{displaymath}
		In that case $Q^{\Gamma}_{\rm KMS}\cong H_{n}^{\infty +}$.
    \eppsn
    {\it Proof}:\\
    If the ${\rm KMS}_{0}$ state is taken to be $\phi(S_{\mu}S_{\nu}^{\ast})=\delta_{\mu,\nu}\frac{1}{n}$, then the two categories coincide by Theorem \ref{newold}. So in that case to determine $Q^{\Gamma}_{\rm KMS}$, it suffices to compute the universal object in the category $\clc_{\tau}^{\rm Lin}$. To that end we denote the universal object by $Q^{\rm Lin}_{\tau}$ and	let $\alpha:C^{\ast}(\Gamma)\raro C^{\ast}(\Gamma)\ot Q^{\rm Lin}_{\tau}$ be the action of $Q^{\rm Lin}_{\tau}$ given by $\alpha(S_{i})=\sum_{j=1}^{n}S_{j}\ot q_{ji}$ for $i=1,...,n$. Then since $\alpha$ preserves $\tau$, it follows easily that both the matrices $((q_{ij}))$ and $((q_{ij}^{\ast}))$ are unitary. Also using the fact that $S_{i}$'s are partial isometries and $S_{i}^{\ast}S_{j}=S_{i}S_{j}^{\ast}=0$ for $i\neq j$,  we can show that $q_{ij}$ is partial isometry  for all $i,j=1,...,n$. Using $S_{i}S_{j}^{\ast}=\delta_{ij}p_{i}=S_{i}^{\ast}S_{j}$ for all $i,j=1,...,n$, we have the following for all $i=1,...,n$:
    \begin{eqnarray*}
    	&&\alpha(S_{i}S_{i}^{\ast})=\alpha(S_{i}^{\ast}S_{i})\\
    	&\Rightarrow& \sum_{k=1}^{n} p_{k}\ot q_{ki}q_{ki}^{\ast}=\sum_{k=1}^{n}p_{k}\ot q_{ki}^{\ast}q_{ki}\\
    	&\Rightarrow&\sum_{k=1}^{n}p_{k}\ot(q_{ki}q_{ki}^{\ast}-q_{ki}^{\ast}q_{ki})=0.
    \end{eqnarray*} 
    Now since $p_{k}$'s are mutually orthogonal projections, it follows that $q_{ki}^{\ast}q_{ki}=q_{ki}q_{ki}^{\ast}$ for all $i,k=1,...,n$. So by definition of $H_{n}^{\infty +}$, there is a surjective CQG morphism from $Q^{\rm Lin}_{\tau}$ to $H_{n}^{\infty +}$ sending generators to generators.\\
    \indent Now we show that $H_{n}^{\infty +}$ is an object in the category $\clc^{\rm Lin}_{\tau}$, so that there is a surjective CQG morphism from $H_{n}^{\infty +}$ to $Q^{\rm Lin}_{\tau}$ sending generators to generators and hence establishing that $Q^{\rm Lin}_{\tau}\cong H^{\infty +}_{n}$. For that define $\alpha(S_{i})=\sum_{j=1}^{n}S_{j}\ot q_{ji}$ for all $i=1,...,n$, $q_{ij}\in H_{n}^{\infty +}$. Note that it suffices to prove that $\alpha$ is a well defined $C^{\ast}$-homomorphism. Rest of the conditions of action of CQG follow immediately. Also since $((q_{ij}^{\ast}))$ is unitary, it follows that $\alpha$ preserves $\tau$. To show that $\alpha$ is a well defined $C^{\ast}$-homomorphism, we use universality of $C^{\ast}(\Gamma)$. Since, $q_{ij}$'s are partial isometries, it follows that $\alpha(S_{i})$'s are partial isometries for all $i=1,...,n$. Also using the fact that $q_{ij}$'s are normal elements with both $q_{ij}q_{ij}^{\ast}$ and $q_{ij}^{\ast}q_{ij}$'s projections, we can show that $\alpha(S_{i})^{\ast}\alpha(S_{i})=\alpha(S_{i})\alpha(S_{i})^{\ast}$ and each of them are self adjoint say $P_{i}$. Hence the proof of the fact that $Q^{\rm Lin}_{\tau}\cong H^{\infty +}_{n}$ will be complete if we show that $P_{i}P_{j}=\delta_{ij}P_{i}$ for all $i,j=1,...,n$.
    \begin{eqnarray*}
    	P_{i}P_{j}&=& \alpha(S_{i})^{\ast}\alpha(S_{i})\alpha(S_{j})\alpha(S_{j})^{\ast}\\
    	&=&(\sum_{k=1}^{n} p_{k}\ot q_{ki}^{\ast}q_{ki})(\sum_{l=1}^{n} p_{l}\ot q_{lj}q_{lj}^{\ast})\\
    	&=&\sum_{k=1}^{n} p_{k}\ot q_{ki}^{\ast}q_{ki}q_{kj}q_{kj}^{\ast}.
    \end{eqnarray*}
    For $i=j$, using the fact that $q_{ki}$ is normal, we see that $P_{i}^{2}=P_{i}$. For $i\neq j$, since by Lemma \ref{H}, $q_{ki}q_{kj}=0$ for all $k=1,..,n$, we have $P_{i}P_{j}=0$. Hence by Corollary \ref{old}, $Q^{\Gamma}_{\rm KMS}\cong H_{n}^{\infty +}$.\\
		\indent For the converse, let us take any ${\rm KMS}_{0}$ state $\psi$ given by any probability measure $(\tau_{1},...,\tau_{n})$ on the vertex set so that $\tau_{i}=\psi(S_{i}S_{i}^{\ast})$ for all $i=1,...,n$. Now we have already seen that $H_{n}^{\infty +}$ is the universal object in the category $\clc^{\rm Lin}_{\tau}$. We shall show that if $H_{n}^{\infty +}$ belongs to the category $\clc^{\Gamma}_{\psi}$, then $\tau_{i}=\tau_{j}$ for all $i,j=1,...,n$, completing the proof of the Proposition. Let $\alpha:C^{\ast}(\Gamma)\raro C^{\ast}(\Gamma)\ot H_{n}^{\infty +}$ be the action given by $\alpha(S_{i})=\sum_{j=1}^{n}S_{j}\ot q_{ji}$. If $H_{n}^{\infty +}$ belongs to the category $\clc^{\Gamma}_{\psi}$, $\alpha$ preserves $\psi$. Hence we have
		\begin{eqnarray*}
		&&(\psi\ot{\rm id})\circ\alpha(S_{i}S_{i}^{\ast})=\psi(S_{i}S_{i}^{\ast}).1\\
		&\Rightarrow& \sum_{k=1}^{n}\tau_{k}q_{ki}q_{ki}^{\ast}=\tau_{i}(\sum_{k=1}^{n}q_{ki}q_{ki}^{\ast})\\
		&\Rightarrow& \sum_{k=1}^{n}(\tau_{k}-\tau_{i})q_{ki}q_{ki}^{\ast}=0.
		\end{eqnarray*}
		Now using the fact that in $H_{n}^{\infty +}$, $q_{ki}q_{ki}^{\ast}$ are orthogonal projections for $k=1,...,n$, we conclude that $\tau_{i}=\tau_{k}$ for all $k=1,...,n$, completing the proof of the Proposition.\qed
    \subsection{Quantum symmetry of Cuntz algebra at critical inverse temperature}
    Recall the definition of Cuntz algebra (to be denoted by $\clo_{n}$) from \cite{ped}. $\clo_{n}$ is the universal $C^{\ast}$-algebra generated by $n$ partial isometries $\{S_{i}:i=1,...,n\}$ satisfying the following relations:
    \begin{eqnarray*}
    	&& S_{i}^{\ast}S_{i}=1, i=1,...,n\\
    	&&\sum_{j=1}^{n}S_{j}S_{j}^{\ast}=1.
    \end{eqnarray*}
    It is a graph $C^{\ast}$ algebra with the vertex matrix a scalar matrix $D=n$ and $\rho(D)=n$. So it is a strongly connected graph which  satisfies our hypothesis on the graph so that it has a unique ${\rm KMS}_{{\rm log} n}$ state $\phi$. It is given on the dense $\ast$-algebra $\cla_{0}$ by
    \begin{displaymath}
    \phi(S_{\mu}S_{\nu}^{\ast})=\delta_{\mu,\nu}\frac{1}{n^{|\mu|}}.
    \end{displaymath}
    This gauge invariant state is very well studied and reader might see \cite{ped} for details on this state. It has an alternative description given by $\phi=\tau\circ\Phi$, where $\Phi$ is the canonical expectation from $\clo_{n}$ to $\clo_{n}^{\gamma}$ and $\tau$ is the tracial state on the AF algebra $\clo_{n}^{\gamma}$. The state $\phi$ is also faithful on $\clo_{n}$.
    \blmma
    \label{plug}
    For $|\mu|=|\nu|$, $\phi(S_{\mu}aS_{\nu}^{\ast})=\delta_{\mu,\nu}\frac{1}{n^{|\mu|}}\phi(a) \ \forall a\in\clo_{n}$.
    \elmma
    {\it Proof}:\\
    Let $a=S_{\delta} S_{\gamma}^{\ast}$. If $\mu\neq\nu$, by definition of $\phi$, $\phi(S_{\mu}S_{\delta}S_{\gamma}^{\ast}S_{\nu}^{\ast})=0$. Otherwise,  if $|\delta|\neq|\gamma|$, by definition of $\phi$, $\phi(S_{\mu}S_{\delta}S_{\gamma}^{\ast}S_{\mu}^{\ast})=0=\phi(S_{\delta}S_{\gamma}^{\ast})$ as $|\mu|=|\nu|$. So let $|\delta|=|\gamma|$. If $\delta\neq\gamma$, $\phi(S_{\mu}S_{\delta}S_{\gamma}^{\ast}S_{\mu}^{\ast})=0=\phi(S_{\delta}S_{\gamma}^{\ast})$. If $\delta=\gamma$, then $\phi(S_{\mu}S_{\delta}S_{\delta}^{\ast}S_{\mu}^{\ast})=\frac{1}{n^{(|\mu|+|\delta|)}}=\frac{1}{n^{|\mu|}}\phi(S_{\delta}S_{\delta}^{\ast})$. Hence for $a=S_{\delta}S_{\gamma}^{\ast}$,  
    \begin{displaymath}
    \phi(S_{\mu}aS_{\nu}^{\ast})=\delta_{\mu,\nu}\frac{1}{n^{|\mu|}}\phi(a).
    \end{displaymath}
    Now by linearity and continuity of $\phi$, we complete the proof of the lemma.\qed\vspace{0.1in}\\
    We have an  innerproduct on $\clo_{n}$ induced by the state $\phi$ given by $\langle a,b\rangle:=\phi(a^{\ast}b)$. We are going to introduce an orthogonal filtration on $\clo_{n}$ with respect to the state $\phi$. To that end let $\mathcal{B}_{r,s}:={\rm Sp}\{S_{\mu}S_{\nu}^{\ast}:|\mu|=r,|\nu|=s\}$. We have the following
    \blmma
    \label{ortho_init}
    For $a\in\mathcal{B}_{r,s}$, $b\in\mathcal{B}_{r^{\prime},s^{\prime}}$, $a\perp b$ if $(r+s^{\prime})\neq(r^{\prime}+s)$.
    \elmma
    {\it Proof}:\\
    Let $r> r^{\prime}$ and $a=S_{\mu_{1}}...S_{\mu_{r}}S_{\gamma_{1}}^{\ast}...S_{\gamma_{s}}^{\ast}\in\mathcal{B}_{r,s}$, $S_{\delta_{1}}...S_{\delta_{r^{\prime}}}S_{\xi_{1}}^{\ast}...S_{\xi_{s^{\prime}}}^{\ast}\in\mathcal{B}_{r^{\prime},s^{\prime}}$. Then
    \begin{displaymath}
    a^{\ast}b=S_{\gamma_{s}}...S_{\gamma_{1}}S_{\mu_{r}}^{\ast}...S_{\mu_{1}}^{\ast}S_{\delta_{1}}...S_{\delta_{r^{\prime}}}S_{\xi_{1}}^{\ast}...S_{\xi_{s^{\prime}}}^{\ast}.
    \end{displaymath}
    If $(\mu_{1},...,\mu_{r^{\prime}})\neq (\delta_{1},...,\delta_{r^{\prime}})$, then $a^{\ast}b=0$. So let us assume that $(\mu_{1},...,\mu_{r^{\prime}})=(\delta_{1},...,\delta_{r^{\prime}})$. Then 
    \begin{displaymath}
    a^{\ast}b=S_{\gamma_{s}}...S_{\gamma_{1}}S_{\mu_{r}}^{\ast}...S_{\mu_{r^{\prime}+1}}^{\ast}S_{\xi_{1}}^{\ast}...S_{\xi_{s^{\prime}}}^{\ast},
    \end{displaymath}
    i.e. $a^{\ast}b\in\mathcal{B}_{s,(r-r^{\prime}+s^{\prime})}$. Hence by definition of $\phi$, $\phi(a^{\ast}b)=0$ if $s\neq (r-r^{\prime}+s^{\prime})$ i.e. $a\perp b$ if $(r+s^{\prime})\neq(r^{\prime}+s)$.\\
    For $r\leq r^{\prime}$, $a^{\ast}b=0$ if $(\mu_{1},...,\mu_{r})\neq(\delta_{1},...,\delta_{r})$. So let us assume that $(\mu_{1},...,\mu_{r})=(\delta_{1},...,\delta_{r})$. Then
    \begin{displaymath}
    a^{\ast}b=S_{\gamma_{s}}...S_{\gamma_{1}}S_{\delta_{r+1}}...S_{\delta_{r^{\prime}}}S_{\xi_{1}}^{\ast}...S_{\xi_{s^{\prime}}}.
    \end{displaymath}
    So $a^{\ast}b\in\mathcal{B}_{(s+r^{\prime}-r),s^{\prime}}$. So by definition of $\phi$, $\phi(a^{\ast}b)=0$ if $(s+r^{\prime}-r)\neq s^{\prime}$ i.e. $a\perp b$ if $(r+s^{\prime})\neq(r^{\prime}+s)$.\qed\vspace{0.1in}\\
    Recall the fixed point subalgebra $\clo_{n}^{\gamma}$ of $\clo_{n}$ with respect to the canonical gauge action $\gamma$ on $\clo_{n}$. $\clo_{n}^{\gamma}$ is a UHF algebra with the finite dimensional $C^{\ast}$-subalgebras  given by 
    \begin{displaymath}
    \clf_{0}=\mathbb{C}1, \clf_{k}:={\rm Sp}\{S_{\mu}S_{\nu}^{\ast}:|\mu|=|\nu|=k\}, k\in\mathbb{N}.
    \end{displaymath}
    Observe that in the notation introduced in the beginning of this Subsection, $\clf_{k}$ is nothing but $\mathcal{B}_{k,k}$. As graph of Cuntz algebra is without sink, recall from Subsection \ref{KMS} that $\clf_{k}\subset\clf_{k+1}$ for all $k\in\mathbb{N}_{0}$. Hence,  we can introduce the following  sequence (see \cite{Evans}) of orthogonal finite dimensional vector subspaces: 
    \begin{displaymath}
    \clw_{0}=\clf_{0}=\mathbb{C}1, \clw_{1}=\clf_{1}\ominus\clf_{0},...,\clw_{k}=\clf_{k}\ominus\clf_{k-1},...
    \end{displaymath}
    Using $\clw_{k}$, we define  finite dimensional vector subspaces for $(k,r)\in \mathbb{N}_{0}\times\mathbb{N}$, where $\mathbb{N}_{0}$ stands for the set $\mathbb{N}\cup \{0\}$.
    \begin{eqnarray*}
    	\clv^{(1)}_{k,r}:={\rm Sp}\{S_{\mu}x:|\mu|=r,x\in\clw_{k}\}\\
    	\clv^{(2)}_{k,r}:={\rm Sp}\{yS_{\gamma}^{\ast}:|\gamma|=r,y\in\clw_{k}\}.
    \end{eqnarray*}
    \bppsn
    $(\phi,(\clw_{k},\clv^{(1)}_{l,r},\clv^{(2)}_{p,s}: k\in\mathbb{N}_{0}, (l,r)\in\mathbb{N}_{0}\times \mathbb{N},(p,s)\in\mathbb{N}_{0}\times\mathbb{N}))$ is an orthogonal filtration for $\clo_{n}$.
    \eppsn
    {\it Proof}:\\
    Recall definition \ref{filtration} of orthogonal filtration. Condition (i) is satisfied, since $\clw_{0}=\mathbb{C}1$. For the span density condition let $S_{\mu}S_{\nu}^{\ast}\in \cla_{0}$. If $|\mu|=|\nu|$, then $S_{\mu}S_{\nu}^{\ast}\in\clf_{|\mu|}$ i.e. $S_{\mu}S_{\nu}^{\ast}\in {\rm Sp}\cup\clw_{k}$. So let $|\mu|>|\nu|$. Then $S_{\mu}S_{\nu}^{\ast}=S_{\delta}S_{\mu^{\prime}}S_{\nu}^{\ast}$ where $|\mu^{\prime}|=|\nu|$. As before $S_{\mu^{\prime}}S_{\nu}^{\ast}\in {\rm Sp}\cup\clw_{k}$ i.e. $S_{\mu}S_{\nu}^{\ast}\in{\rm Sp}\cup\clv^{(1)}_{k,r}$ if $|\delta|=r>0$. Similarly we can show that if $|\mu|<|\nu|$, $S_{\mu}S_{\nu}^{\ast}\in {\rm Sp}\cup\clv^{(2)}_{k,r}$ for some $(k,r)\in\mathbb{N}_{0}\times \mathbb{N}$. Using the density of $\cla_{0}$ in $C^{\ast}(\Gamma)$, density of ${\rm Sp}(\cup\clw_{k}\cup\clv^{(1)}_{l,r}\cup\clv^{(2)}_{p,s})$ in $C^{\ast}(\Gamma)$ is clear now. Now we proceed to prove the second condition of definition \ref{filtration}. We start with the simple observation that $\clw_{k}\subset \clf_{k}=\mathcal{B}_{k,k}$, $\clv^{(1)}_{k,r}\subset \clb_{r+k,k}$ and $\clv^{(2)}_{k,r}\subset\clb_{k,r+k}$. By construction $\clw_{k}\perp\clw_{l}$ for $k\neq l$. For $a\in \clw_{k}$, $b\in\clv^{(1)}_{l,r}$, $\phi(a^{\ast}b)=0$ by Lemma \ref{ortho_init}, since $r\in\mathbb{N}$ implies $(k+l+r)\neq (k+l)$. With similar reasoning, for $a\in\clw_{k}$ and $b\in\clv^{(2)}_{l,r}$, $a\perp b$. Now let $a\in\clv^{(1)}_{k,r}, b\in\clv^{(2)}_{l,s}$. So $a\in \clb_{r+k,k}$ and $b\in\clb_{l,l+s}$. So another application of Lemma \ref{ortho_init} shows that $\phi(a^{\ast}b)=0$, since $(r+k+l+s)=(k+l)$ implies $r=-s$ which is not possible since $(r,s)\in\mathbb{N}\times\mathbb{N}$. Now let $a\in\clv^{(1)}_{k,r}, b\in\clv^{(1)}_{l,s}$ for $(k,r)\neq(l,s)$. If $r\neq s$, $\phi(a^{\ast}b)=0$ by definition of $\phi$. If $r=s$, $k\neq l$, let $a=S_{\mu}x\in \clv^{(1)}_{k,r}$ for $|\mu|=r,x\in\clw_{k}$ and $b=S_{\nu}y\in\clv^{(1)}_{l,r}$ for $|\nu|=r, y\in\clw_{l}$. Then $\phi(a^{\ast}b)=\delta_{\mu,\nu}\phi(x^{\ast}y)=0$, since $k\neq l$. Similarly for $a=xS_{\mu}^{\ast}\in\clv^{(2)}_{k,r},b=yS_{\nu}^{\ast}\in\clv^{(2)}_{l,s}$, $\phi(a^{\ast}b)=\phi(S_{\mu}x^{\ast}yS_{\nu}^{\ast})$. If $r\neq s$, $\phi(a^{\ast}b)=0$. If $r=s$, by Lemma \ref{plug},  $\phi(S_{\mu}x^{\ast}yS_{\nu}^{\ast})=0$ since $\phi(x^{\ast}y)=0$ for $k\neq l$. This completes the proof of the proposition.\qed\vspace{0.1in}\\
    Let $\clc_{\clo_{n},\mathcal{V}}$ be the category whose objects are $(\alpha,(\clq,\Delta))$, $\alpha$ is an action of $(\clq,\Delta)$ such that it preserves the orthogonal filtration $(\phi,(\clw_{k},\clv^{(1)}_{l,r},\clv^{(2)}_{p,s}: k\in\mathbb{N}_{0}, (l,r)\in\mathbb{N}_{0}\times \mathbb{N},(p,s)\in\mathbb{N}_{0}\times\mathbb{N}))$ on $\clo_{n}$. Then we have the following
    \bthm
    If we denote the graph of $\clo_{n}$ by  $\Gamma$, then the categories $\clc^{\Gamma}_{\phi}$ and $\clc_{\clo_{n},\mathcal{V}}$ are isomorphic.
    \ethm
    {\it Proof}:\\
    Let $(\alpha,(\clq,\Delta))$ be an object in the category $\clc_{\clo_{n},\clv}$. Then $\alpha$ preserves the state $\phi$ by remark \ref{statepreserve} in the Subsection \ref{orthogonalfiltration}. Also note that $\clv^{(1)}_{0,1}={\rm Sp}\{S_{i}:i=1,...,n\}$ which implies that $\alpha$ is linear. Hence $(\alpha,(\clq,\Delta))$ belongs to the category $\clc^{\Gamma}_{\phi}$.\\
    \indent Conversely suppose $(\alpha,(\clq,\Delta))$ is an object in the category $\clc^{\Gamma}_{\phi}$. We shall show that $\alpha$ preserves the orthogonal filtration $(\phi,(\clw_{k},\clv^{(1)}_{l,r},\clv^{(2)}_{p,s}: k\in\mathbb{N}_{0}, (l,r)\in\mathbb{N}_{0}\times \mathbb{N},(p,s)\in\mathbb{N}_{0}\times\mathbb{N}))$.\vspace{0.1in}\\
    {\bf Step 1}: In step 1, we shall show that $\alpha(\clw_{k})\subset\clw_{k}\ot\clq$ for all $k$. Let $\alpha$ be given by $\alpha(S_{i})=\sum_{j=1}^{n}S_{j}\ot q_{ji}$ for $i=1,...,n$. Since $\alpha$ is linear, $\alpha(\clf_{k})\subset\clf_{k}\ot\clq$ for all $k$ i.e. $\alpha$ is a (co)representation of the CQG $\clq$ on the finite dimensional vector space $\clf_{k}$ for all $k$. Also $\alpha$ preserves $\phi$ which, restricted to $\clf_{k}$, is nothing but the matrix trace $\tau$. Hence $\alpha$ is actually a unitary representation on $\clf_{k}$. Since it preserves the finite dimensional subspace $\clf_{k-1}$ of $\clf_{k}$, it also preserves $\clf_{k}\ominus\clf_{k-1}$ which is nothing but $\clw_{k}$. Hence $\alpha(\clw_{k})\subset \clw_{k}\ot\clq$ for all $k$.\vspace{0.1in}\\ 
    {\bf Step 2}: In step 2, we shall show that $\alpha(\clv^{(i)}_{k,r})\subset\clv^{(i)}_{k,r}\ot\clq$, for $i=1,2$ and $(k,r)\in\mathbb{N}_{0}\times\mathbb{N}$, completing the proof of the theorem. To that end let $S_{\mu_{1}}...S_{\mu_{r}}x\in \clv^{(1)}_{k,r}$ where $x\in\clw_{k}$. Hence $\alpha(S_{\mu_{1}}...S_{\mu_{r}}x)=(\sum S_{i_{1}}...S_{i_{r}}\ot q_{i_{1}\mu_{1}}...q_{i_{r}\mu_{r}})(\alpha(x))$. Now, by step 1, $\alpha(x)\in\clw_{k}\ot\clq$ proving that $\alpha(S_{\mu_{1}}...S_{\mu_{r}}x)\in\clv^{(1)}_{k,r}\ot\clq$ i.e. $\alpha(\clv^{(1)}_{k,r})\subset \clv^{(1)}_{k,r}\ot\clq$. Following exactly similar argument, it is easy to see that $\alpha(\clv^{(2)}_{k,r})\subset\clv^{(2)}_{k,r}\ot\clq$, completing the proof of the theorem.\qed\vspace{0.1in}\\
    Combining the above theorem with Corollary \ref{old} and Proposition 5.9 of \cite{sou_arn}, we have
    \bcrlre
    If we denote the universal object  in the category $\clc_{\clo_{n},\clv}$ by $Q_{\clo_{n},\clv}$, then $Q_{\clo_{n},\clv}\cong Q^{\Gamma}_{\rm KMS}\cong Q^{\rm Lin}_{\tau}\cong U^{+}_{n}$.
    \ecrlre
    {\bf Acknowledgement}: We would like to thank Debashish Goswami for encouraging us to look into this problem as well as for some fruitful discussions. We also thank Adam Skalski for clarifying a few points regarding orthogonal filtration.

Soumalya Joardar \\
Theoretical Science Unit,\\ 
JNCASR, Bangalore-560064, India\\ 
email: soumalya.j@gmail.com 
\vspace{0.1in}\\
Arnab Mandal\\
School Of Mathematical Sciences\\
NISER, HBNI,  Bhubaneswar,  Jatni-752050, India\\
email: arnabmaths@gmail.com

\end{document}